# Behavior of power plants in collusion and competitive markets considering external costs

Mehdi Peyrovan[1], Sadoullah Ebrahimnejad[2, *], Amirhossein Moosavi[2]

[1]*Department of Industrial Engineering, South-Tehran Branch, Islamic Azad University, Tehran, Iran, m.pey.teach@gmail.com.*

[2, *]*Department of Industrial Engineering, Karaj Branch, Islamic Azad University, Karaj, Iran, ibrahimnejad@kiau.ac.ir; saeed_ebrahimnejad@yahoo.com,* Tel: (+98) 912-2099648.

[2] *Telfer School of Management, University of Ottawa, 55 Laurier Avenue East, Ottawa, Ontario K1N 6N5, Canada, amir.moosavi@uottawa.ca.*

***Corresponding author***: Sadoullah Ebrahimnejad

Postal Address (corresponding author): Department of Industrial Engineering, Islamic Azad University, Karaj Branch, Rajaee Shahr, Karaj, Iran. P.O. Box: 31485/313. Telephone number: +98 912 2099648.

# Behavior of power plants in collusion and competitive markets considering external costs

**Abstract:** In this paper, a mathematical model is proposed to optimize the production level of Power Plants (PPs) in collusion and competitive markets and in the presence of external costs. To evaluate the model, two meta-heuristics, such as Genetic Algorithm (GA), and Particle Swarm Optimization (PSO), are proposed. A numerical example is studied to assess the behavior of PPs in the presence of external costs. According to the computational results, the increment in the external costs will decrease the production level and profit of all PPs in both collusion and competitive markets. But the computational results show that the effect of increment in the external costs on the PPs with high-quality equipment is less compared to the PPs with poor equipment. Moreover, the increment in the external costs will force PPs to use greener fuel resources instead of fossil fuel.

**Keywords:** Power Plant; Environmental Pollutant; External Cost; Collusion Market; Competitive Market; Genetic Algorithm; Particle Swarm Optimization.

## 1. Introduction

In today's world, power industry plays a fundamental role in daily life of people and provides necessary facilities for development of countries. Since Michael Faraday introduced the underlying principles of electricity generation in the 1820s, PPs use the same principles so far, i.e., generating electricity by converting other resources.

Because PPs use fossil fuel as their main resources to generate electricity, they would probably emit a significant amount of environmental pollutant. Due to increasing consciousness of people in all around the world, governments are enacting some regulations and providing incentive plans for PPs to protect the environment. Thus, decision-makers in this industry have been motivated to control and decrease the environmental contaminants emitted by PPs. The process of converting primary resources to electricity emits hazardous pollutants, such as $CO_2$, $CO$, $SO_2$, $SO_3$, and $NO_x$. Emission of such environmental contaminants imposes noticeable harms to the environment. Therefore, this problem has drawn the attention of researchers in recent years (Vrhovcak et al., 2005). In contrast to the most of developed countries, the majority of under-developed countries have not paid attention to hazardous effects of environmental pollutants emitted by PPs into the environment. In such a situation, the concurrent price of electricity in markets is nominal and does not include the external costs that are imposed on the environment (Weinzettel et al., 2012). Adding these external costs to the nominal price of electricity would be a necessary and effective step toward motivating PPs to take responsibility for the environmental pollutants. In addition, considering external costs could provide the conditions required for supplying sustainable and greener energies (Khoshakhlagh et al., 2013).

In the power industry, competition between PPs might be risky due to high-level of investment (Jaehnert and Doorman, 2012). Considering the fact that there are only a few number of PPs in the majority of under-developed countries, and moreover, homogenous nature of electricity, electricity market can be called as collusion or competitive markets. In such markets, PPs can pursue various policies, e.g., decrease their production level, in order to increase the price of electricity (Stoft, 2002). Thus, cooperation mode with other companies is fundamental in such markets. In addition, it is worth mentioning that this decision affects market share and profit acquired by companies. Hence, the behavior of PPs in collusion and competitive markets has been investigated and compared in this paper. For this purpose, a non-linear optimization model for PPs, which considers abovementioned external costs, has been proposed in this study. The proposed mathematical model aims to determine the production level of PPs in the presence of external costs. To the best of our knowledge, merely few studies have considered external costs for PPs.

In this study, two meta-heuristic algorithms, such as GA and PSO, have been proposed to solve large-scale instances of the problem. Then, the performance of these two algorithms have been compared with each other to determine which of them are more efficient in terms of solution quality and CPU Time.

The remainder of this paper is organized as follows; following the literature review in Section 2, mathematical models for collusion and competitive markets are proposed in Section 3. The meta-heuristic algorithms are developed in Section 4. Then, a numerical example is studied in Section 5. Eventually, Section 6 presents the conclusion.

## 2. Literature review

As previously mentioned, colluding or competing with each other is a critical decision for PPs. Carried our investigation by Sweeting (2007) showed that Britain PPs has colluded with each other between 1992 and 2000. Fabra and Toro (2005) also showed that PPs in the decentralized market of Spain had colluded with each other in 1998. In such a market, PPs have to consider the mutual interests of each other and cooperate with each other. Therefore, they try to maximize their shared profits with respect to the market's demand. This kind of cooperation, in which consumers confronts only with one electricity generator, is called Cartel. Jehle (2001) proposed a mathematical model for collusion market which maximizes all profit functions of PPs.

Unlike the cooperation mode of PPs in collusion market, PPs pursue conflicting goals in the competitive market. Consequently, each PP considers others as a competitor. Given that source of profit is limited in this condition, Nash equilibrium has been used in most of the previous studies in order to find the optimal policy for PPs. In the case of agreeing on one of the possible solutions, PPs gain more profit in comparison to not agreeing on any solution. There are many investigations in the literature relevant to the competitive market for PPs, such as Wen and David (2001), Hobbs (2001), Li and Shahidehpour (2005), Lise et al. (2006), Yucekaya et al. (2009), Li et al. (2011), Fan et al. (2012), Shafie-Khah et al. (2013), Ladjici et al. (2014).

Cournot, Bertrand, and Supply function equilibrium are relevant models to collusion and competitive markets. According to Willems et al. (2009), researchers mostly have used Cournot and Supply function equilibrium in order to determine the price of each unit of electricity. Cournot model considers power market as an energy pool which PPs maximizes their profit by increasing or decreasing the production level of themselves. Based on the supply of PPs, the price of electricity is determined in the Cournot model. In the other word, PPs determine the price of electricity by determining their production level. Some studies, such as Hobbs (2001), Willems (2002), Ellersdorfer (2005), and Bushnell et al. (2007), have used Cournot model in order to determine the price of electricity. On the other hand, PPs maximize their profit based on supply curve (dependent on the price) in the Supply function equilibrium. Unlike the Cournot model, the supply of each PP does not depend on the supply of other PPs in the Supply function equilibrium. Moitre (2002), Evans and Green (2005), Li and Shahidehpour (2005), Holmberg (2008), and Shafie-Khah et al. (2013).

As noted, generating electricity imposes hazardous external costs to the environment, human health, and climate. Respecting to environment (i.e., damage to agriculture products, industrial goods, and species of animals), human health, and climate (i.e., global warming, heat waves, storms, and rising sea levels), Streimikiene and Alisauskaite-Seskiene (2014) evaluated the emission costs of environmental contaminant, i.e., $CO_2$, $SO_2$, and $NO_x$, for PPs in Lithuania. Table (1) shows the determined external costs by Streimikiene and Alisauskaite-Seskiene (2014).

Table 1. External costs of pollutant emission based on Streimikiene and Alisauskaite-Seskiene (2014) - (USD/gr)

| Pollutant | Environment | | | Health | Climate | Total Cost |
| --- | --- | --- | --- | --- | --- | --- |
| | Damage to biodiversity | Damage to product | Damage to materials | | | |
| $NO_x$ | $0.985 \times 10^{-3}$ | $0.357 \times 10^{-3}$ | $0.077 \times 10^{-3}$ | $6.098 \times 10^{-3}$ | - | $7.517 \times 10^{-3}$ |
| $SO_2$ | $0.193 \times 10^{-3}$ | $-0.029 \times 10^{-3}$ | $0.283 \times 10^{-3}$ | $6.621 \times 10^{-3}$ | - | $7.068 \times 10^{-3}$ |
| $CO_2$ | - | - | - | - | $0.023 \times 10^{-3}$ | $0.023 \times 10^{-3}$ |

Considering the high importance of environment, human health, and climate, there are also more investigations in the literature which studied the external costs of pollutants' emission, such as Roth and Ambs (2004), Vrhovcak et al. (2005), Weinzettel et al. (2012), and Georgakellos (2012).

Regarding the abovementioned studies, it can be claimed that there merely few studies that have considered external costs of pollutant emission; furthermore, to best of our knowledge, there are fewer studies that have compared the behavior of PPs in collusion and competitive markets.

## 3. The proposed mathematical model

In this study, we have proposed two mathematical formulations for PPs in both collusion and competitive markets. These mathematical formulations aim to determine the production level of PPs in order to maximize their profit. Besides, these models minimize establishing and maintenance costs of PPs, purchasing fuel resources, and external costs of pollutant emission.

First of all, we have introduced the main principles used for a PP in the following. Then, we have adapted them to multiple PPs. To do this, we have introduced indices, parameters, and decision variables as follows:

- **Indices:**

| | | |
| --- | --- | --- |
| $j$ | Type of fuel resources | $j \in \{1, 2, \dots, J\}$ |
| $k$ | Type of environmental contaminants | $k \in \{1, 2, \dots, K\}$ |

- **Parameters:**

| | |
| --- | --- |
| $\alpha, \beta, \gamma$ | Efficiency coefficients of a PP |
| $H_j$ | Reverse of heating value for $j^{th}$ type of fuel resources |
| $\mu$ | Waste coefficient of a PP |

| $\delta', \delta$ | Coefficients of Cournot model |
| --- | --- |
| $w_j$ | Price of $j^{th}$ type of fuel resources |
| $C_{f\&o\&m}$ | Costs of establishing, generating, and maintenance per each unit of electricity generated in a PP |
| $\vartheta_{jk}$ | Amount of $k^{th}$ type contaminant emitted per each unit of $j^{th}$ type of fuel resources |
| $ec_k$ | External cost of $k^{th}$ type of contaminant |
| $p_{(max)}$ | Maximum capacity of a PP throughout a year |
| $\sigma_j$ | Available amount of $j^{th}$ type of fuel resources throughout a year |
| $Z_k$ | Permissible threshold for emission of $k^{th}$ type of contaminant |
| $s$ | Coefficient of paid subsidies to a PP |

- **Decision variables:**

| $p_j$ | Amount of electricity generated in a PP using $j^{th}$ type of fuel |
| --- | --- |
| $F_j$ | Energy required for a PP in order to generate electricity using $j^{th}$ type of fuel resources |
| $\rho$ | Price of electricity unit |
| $Sub$ | Subside acquired by a PP from government |

In the current investigation, the energy required for PPs to generate electricity by $j^{th}$ type fuel resources has been calculated based on the quadratic function proposed by El-Hawary (1995):

$$F_j = \alpha p_j^2 + \beta p_j + \gamma \qquad (1)$$

To implement the external costs, the Cournot model has been applied to determine the real price of generating each unit of electricity. In this model, the price of each unit of electricity depends on the pure amount of electricity generated by all PPs, i.e., total amount of electricity generated minus waste. Note that the waste of each PP has also been calculated based on the quadratic function proposed by Grainger and Stevenson (1994). The following equation calculates the price of each unit of electricity generated by a PP based on the Cournot model:

$$\rho = \delta - \delta'\left(\sum_{j=1}^{J} p_j - \mu \sum_{j=1}^{J} p_j^2\right) \qquad (2)$$

It has been assumed that a PP acquire subsides from governments based on the following linear equation:

$$Sub = s \cdot \left(\sum_{j=1}^{J} p_j - \mu \sum_{j=1}^{J} p_j^2\right) \qquad (3)$$

Thus, the proposed mathematical model for a PP can be written as follows:

$$Max\ BF = \left(\sum_{j=1}^{J} p_j - \mu \sum_{j=1}^{J} p_j^2\right)\left(\delta - \delta'\left(\sum_{j=1}^{J} p_j - \mu \sum_{j=1}^{J} p_j^2\right)\right) + s \cdot \left(\sum_{j=1}^{J} p_j - \mu \sum_{j=1}^{J} p_j^2\right) -$$
$$\sum_{j=1}^{J} w_j H_j\left(\alpha p_j^2 + \beta p_j + \gamma\right) - \sum_{k=1}^{K} \sum_{j=1}^{J} ec_k \vartheta_{jk} H_j\left(\alpha p_j^2 + \beta p_j + \gamma\right) - C_{f\&o\&m} \sum_{j=1}^{J} p_j \qquad (4)$$

s.t:

$$\sum_{j=1}^{J} p_j \leq p_{(max)} \qquad (5)$$

$$H_j(\alpha p_j^2 + \beta p_j + \gamma) \leq \sigma_j \qquad \forall j \qquad (6)$$

$$\sum_{j=1}^{J} \vartheta_{jk} H_j (\alpha p_j^2 + \beta p_j + \gamma) \leq Z_k \qquad \forall\, k \qquad (7)$$

$$p_j \geq 0 \qquad \forall\, j \qquad (8)$$

The objective function consists of five terms; first and second terms maximize the income and subsides earned by a PP, respectively. The third, fourth, and fifth terms also minimize the costs of purchasing fuel resources required for generation of electricity, pollutants emission, and building, generating and maintenance of a PP, respectively. In order to calculate the average cost of establishing, generating and maintenance of a PP for a specific amount of electricity, total costs of building, generating and maintenance of a PP has to be determined for its whole lifetime considering interest rate. Afterward, prorate this cost to its maximum capacity, i.e., $p_{(max)}$.

Constraint (5) ensures that a PP cannot generate electricity more than its capacity. Constraint (6) guarantees that a PP cannot use more fuel resource of $j^{th}$ type than the predetermined threshold by the government, i.e., $\sigma_j$. Constraint (7) ensures that a PP cannot also emit environmental pollutant more than the predetermined threshold by the government, i.e., $Z_k$. Eventually, Constraint (8) indicates the positive decision variables.

In the following two sub-sections, this model will be adapted for multiple PPs in collusion and competitive markets, respectively.

### 3.1. Adapted model for power plants in collusion market

Considering the homogenous nature of electricity, and a few number of PPs, electricity market can be called as collusion or competitive markets. In this sub-section, we intend to adapt the simple mathematical model proposed earlier for multiple PPs in the collusion market. Note that we have used $i$ as the index of PPs.

$$Max\ (\sum_{i=1}^{I} BF_i) = \sum_{i=1}^{I}(\sum_{j=1}^{M} p_{ij} - \mu_i \sum_{j=1}^{J} p_{ij}^2)\left(\delta - \delta'(\sum_{j=1}^{J} p_{ij} - \mu_i \sum_{j=1}^{J} p_{ij}^2)\right) +$$
$$\sum_{i=1}^{I} s \cdot (\sum_{j=1}^{J} p_{ij} - \mu_i \sum_{j=1}^{J} p_{ij}^2) - \sum_{i=1}^{I}\sum_{j=1}^{J} w_j H_j(\alpha_i p_{ij}^2 + \beta_i p_{ij} + \gamma_i) -$$
$$\sum_{i=1}^{I}\sum_{k=1}^{K}\sum_{j=1}^{J} ec_k \vartheta_{jk} H_j(\alpha_i p_{ij}^2 + \beta_i p_{ij} + \gamma_i) - C_{f\&o\&m} \cdot \sum_{i=1}^{I}\sum_{j=1}^{J} p_{ij} \qquad (9)$$

s.t:

$$\sum_{j=1}^{J} p_{ij} \leq p_{i(max)} \qquad \forall\, i \qquad (10)$$

$$\sum_{i=1}^{I} H_j(\alpha_i p_{ij}^2 + \beta_i p_{ij} + \gamma_i) \leq \sigma_j \qquad \forall\, j \qquad (11)$$

$$\sum_{i=1}^{I}\sum_{j=1}^{J} \vartheta_{jk} H_j(\alpha_i p_{ij}^2 + \beta_i p_{ij} + \gamma_i) \leq Z_k \qquad \forall\, k \qquad (12)$$

$$p_{ij} \geq 0 \qquad \forall\, i, j \qquad (13)$$

It should be noted that Eq. (9) is the collusion or Cartel model for PPs, which tries to maximize the profit earned by all PPs in the collusion market.

### 3.2. Adapted model for power plants in competitive market

In the case of the competitive market, interests of PPs conflict with each other, and they generate electricity separately. In such a circumstances, agreeing on one of the possible cases would be more beneficial for all sides rather than not agreeing upon a case. Therefore, we have used the Nash equilibrium in order to find the best policy for PPs. Using the Nash equilibrium, the proposed model has been adapted for the competitive market, and reformulated as follows:

$$Max\ (BF_1 \times ... \times BF_I) = [(\sum_{j=1}^{J} p_{1j} - \mu_1 \sum_{j=1}^{J} p_{1j}^2)(\delta - \delta'(\sum_{j=1}^{J} p_{1j} - \mu_1 \sum_{j=1}^{J} p_{1j}^2)) +$$

$$s \cdot (\sum_{j=1}^{J} p_{1j} - \mu_1 \sum_{j=1}^{J} p_{1j}^2) - \sum_{j=1}^{J} w_j H_j(\alpha_1 p_{1j}^2 + \beta_1 p_{1j} + \gamma_1) -$$

$$\sum_{k=1}^{K} \sum_{j=1}^{J} ec_k \vartheta_{jk} H_j(\alpha_1 p_{1j}^2 + \beta_1 p_{1j} + \gamma_1) - C_{f\&o\&m} \sum_{j=1}^{J} p_{1j}] \times ... \times$$

$$[(\sum_{j=1}^{j} p_{Nj} - \mu_N \sum_{j=1}^{J} p_{Nj}^2)(\delta - \delta'(\sum_{j=1}^{J} p_{Nj} - \mu_N \sum_{j=1}^{J} p_{Nj}^2)) +$$

$$s \cdot (\sum_{j=1}^{J} p_{Nj} - \mu_N \sum_{j=1}^{J} p_{Nj}^2) - \sum_{j=1}^{J} w_j H_j(\alpha_N p_j^2 + \beta_N p_{Nj} + \gamma_N) -$$

$$\sum_{k=1}^{K} \sum_{j=1}^{J} ec_k \vartheta_{jk} H_j(\alpha_N p_{Nj}^2 + \beta_N p_{Nj} + \gamma_N) - C_{f\&o\&m} \sum_{j=1}^{J} p_{Nj}] \quad (14)$$

s.t:

$$\sum_{j=1}^{J} p_{ij} \leq p_{i(max)} \qquad \forall\ i \quad (15)$$

$$\sum_{i=1}^{I} H_j(\alpha_i p_{ij}^2 + \beta_i p_{ij} + \gamma_i) \leq \sigma_j \qquad \forall\ j \quad (14)$$

$$\sum_{i=1}^{I} \sum_{j=1}^{J} \vartheta_{jk} H_j(\alpha_i p_{ij}^2 + \beta_i p_{ij} + \gamma_i) \leq Z_k \qquad \forall\ k \quad (16)$$

$$p_{ij} \geq 0 \qquad \forall\ i,j \quad (17)$$

## 4. The proposed GA and PSO

The exact methods cannot efficiently solve NP-hard problems (mostly large-scale instances) with respect to CPU Time. Furthermore, the proposed mathematical model is non-linear, and there would be no guarantee for the exact methods to solve it to optimality. Therefore, we have developed two meta-heuristic algorithms, i.e., GA and PSO, to resolve these problems.

### 4.1. Solution representation

One of the most important parts in designing efficient meta-heuristic algorithm is solution representation (the way that solutions are coded). In this study, each solution is represented by **_I_** sections (number of PPs), which each consists of **_J_** subsections ( the type of fuel resources). Random numbers between zero and one are generated and assigned to each of these sub-sections. Then, amount of electricity generated by each type of fuel resources is calculated based on Eq. (18).

$$p_{ij} = (\frac{rand_j}{\sum_{j \in J_i} rand_j}) \times p_{i(max)} \qquad \forall\ i \quad (18)$$

In order to clarify the abovementioned structure, Fig. (1) is provided below for a problem with three PPs and three types of fuel resources. It should be noted that the presented solution representation has been used in both of the proposed meta-heuristic algorithms.

>>Please insert Fig. (1) about here>>

### 4.2. The proposed GA

GA is a recognized population-based meta-heuristic algorithm for computational models inspired by evolution. This algorithm encodes possible solutions of problems on simple structures called chromosome, and apply various operators to these structures in order to create high-quality solutions (Whitley, 1994).

*4.2.1. Initial population and parents selection*

As stated, GA is a population-based meta-heuristic algorithm and require an initial population. For this purpose, we have generated initial population for this algorithm randomly. All chromosomes generated in the first generation are considered as parents, but it is different for the rest of generations. Selecting the parents of the second and subsequent generations is done through tournament method in order to maintain the diversity of the population.

*4.2.2. Fitness function*

After generating the initial population or generating new population in the next iterations, the fitness of each chromosome will be evaluated based on Eqs. (9) and (14) for the collusion or competitive markets, respectively. Note that we have applied the penalty function method in this study to control the infeasible solutions. Further clarification has been provided in sub-section 4.4.

*4.2.3. Crossover and mutation operators*

Crossover and mutation are two important operators in GA. As the most important operator in GA, the main intention of crossover operator is to improve the quality of solutions. Besides, the mutation operator aims to increase the diversity of the population. In this study, two-point crossover and swap mutation have been applied in the proposed GA. Fig. (2) shows two illustrative examples of the crossover and mutation operators.

These operators are applied iteratively to enhance the quality of solutions in the population. This procedure lasts until the termination criterion is met. In this paper, we have defined the termination criterion based on the maximum number of iteration.

>>Please insert Fig. (2) about here>>

## 4.3. The proposed PSO

PSO is also a population-based meta-heuristic algorithm, which was introduced by Kennedy and Eberhart (1995). This algorithm is well-known for its sound performance in dealing with optimization problems. We have developed this algorithm besides GA to validate the performance of both algorithms. Furthermore, we were motivated to figure out which of these algorithms could outperform the other one.

We have originated the settings used in the proposed PSO from Clerc and Kennedy (2002). To sake for brevity, therefore, we refer interested readers to Clerc and Kennedy (2002) for more details.

## 4.4. Penalty function

Given that the proposed meta-heuristic algorithms generate the initial population and apply operators randomly, it is possible that some solutions would be infeasible, i.e., solutions cannot satisfy the constraints of the problem. For this purpose, we have introduced two penalty functions to calculate the violations of solutions.

$$V1 = (\frac{(ep_k - Z_k + 1)^+}{ep_k - Z_k + 1} \times \frac{ep_k}{Z_k - 1}) \times 10^5 \qquad \forall\, k \qquad (19)$$

$$V2 = (\frac{(cr_j - \sigma_j + 1)^+}{cr_j - \sigma_j + 1} \times \frac{cr_j}{\sigma_j - 1}) \times 10^5 \qquad \forall j \qquad (20)$$

where $ep_k$ and $cr_j$ refer to the emission of $k^{th}$ type of pollutants, and consumption of $j^{th}$ type of fuel resources, respectively. Accordingly, Eqs. (19) and (20) calculate penalty costs for emission of environmental contaminants and consumption of fuel resources, respectively. These values are multiplied by an adequately large coefficient, then, added to the main objective functions, i.e., Eqs. (9) and (14).

## 5. Numerical example

In this section, a numerical example has been studied. Authors have investigated the power industry in Iran and has created the numerical example contemporary condition in this country. This numerical example consists of three different PPs. Given that PPs of Iran mostly use fuel-oil, gas-oil, and gas as their primary fuel resources, therefore, PPs would mainly emit $CO_2$, $SO_2$, and $NO_x$ as three hazardous environmental contaminants. Thus, external costs of three mentioned contaminants have been considered in this study.

According to Six-Bus system, we have extracted the efficiency coefficients for all PPs. Afterward, we have specified the maximum capacity of PPs based on their efficiency coefficients. Table (2) shows efficiency coefficients and maximum capacity of PPs.

Table 2. Efficiency coefficients and maximum capacity of PPs

| PP | $\alpha$ | $\beta$ | $\gamma$ | Maximum capacity ($p_{i(\max)}$) |
|---|---|---|---|---|
| 1 | 0.00041 | 15.5 | 1078 | $2.75 \times 10^6$ (MWh) |
| 2 | 0.00031 | 16 | 14 | $2.75 \times 10^6$ (MWh) |
| 3 | 0.00051 | 14 | 702.9 | $2.75 \times 10^6$ (MWh) |

As previously mentioned, PPs mostly use fuel-oil, gas-oil, and gas as their primary fuel resources in Iran. Based on the balance sheet of Iran's energy in 2011, the price of these fuel resources has been determined. Moreover, reverse of the heating value for fuel resources has been extracted from the detailed stats of Iran's electricity industry in 2011. Table (3) displays the price and reverse of the heating value for fuel resources.

Table 3. Price and reverse of the heating value for different types of fuel resources

| Fuel | Fuel-oil (j=1) | Gas-oil (j=2) | Gas (j=3) |
|---|---|---|---|
| Price ($w_j$) | 0.057 (USD/lit) | 0.1 (USD/lit) | 0.022 (USD/m³) |
| Reverse of heating value ($H_j$) | 0.108 (lit/Mcal) | 0.116 (lit/Mcal) | 0.114 (m³/Mcal) |

With respect to the heating value, pollution level of fuel resources has also been determined and provided in Table (4).

Table 4. Pollution level of fuel resources

| Pollution level / Fuel resources | $NO_x$ ($k=1$) | $SO_2$ ($k=2$) | $CO_2$ ($k=3$) |
|---|---|---|---|
| Fuel-oil ($\vartheta_{1k}$) | 5 (gr/lit) | 46.9 (gr/lit) | 2978 (gr/lit) |
| Gas-oil ($\vartheta_{2k}$) | 5.2 (gr/lit) | 15.7 (gr/lit) | 2648 (gr/lit) |
| Gas ($\vartheta_{3k}$) | 3.1 (gr/m³) | 0 (gr/m³) | 2133 (gr/m³) |

In order to calculate the average cost of establishing, generating and maintenance of a PP for a specific amount of electricity, total costs of building, generating and maintenance of a PP has to be determined for its whole lifetime considering interest rate. Afterward, prorate this cost to its maximum capacity. We have originated this cost from other studies in the literature. According to Lajevardi and Mohades (2010), we have set $C_{f\&o\&m}$ equal to 7.1×10⁻³ USD per generating each unit of electricity. It is worth mentioning that each unit of electricity in this study is set equal to Kilowatt hour (KWh).

$\delta$ is a constant coefficient in the Cournot model that always considered to be larger than the cost of generating each unit of electricity. According to the balance sheet of Iran's energy in 2011, cost of generating each unit of electricity has been estimated to be equal to 0.035 USD in Iran; therefore, $\delta$ has been set to be equal to 0.039 USD. $\delta'$ also is one another constant coefficient in the Cournot model. While Chen et al. (2006) and Cunningham et al. (2002) considered this coefficient equal to $10^{-2}$, we also have set $\delta'$ equal to $10^{-2}$ in this study.

Considering the fact that a portion of electricity generated by PPs might be wasted, entire electricity generated by PPs would not be transferred to the market. Amount of waste directly depends on the distance of PPs from the demand markets. Furthermore, the waste coefficients of electricity generated by PPs could vary from each other. To sake for simplicity, we have set the waste coefficient of all three PPs equal to $10^{-8}$.

In the literature, many investigations have evaluated the costs of environmental contaminants. Roth and Ambs (2004) studied destructive effects of environmental pollutants, such as $CO_2$, $SO_2$, and $NO_x$, and estimated three cost ranges for their emission, i.e., lower range, best estimate range, upper range. As they stated, the external costs for emission of $CO_2$, $SO_2$, and $NO_x$ can utmost be equal to $7.1 \times 10^{-3}$, $3.5 \times 10^{-3}$, and $0.028 \times 10^{-3}$ USD per gram, respectively. To assess how gradually implementation of the external costs by the government would affect the market, six different scenarios have been introduced with respect to the upper range of external costs mentioned earlier. These scenarios have been provided in Table (5). Moreover, available fuel resources for generating electricity, and permissible amount for emission of pollutants are introduced in Tables (6) and (7), respectively.

Table 5. Six scenarios of implementing the external costs by the government - (USD/gr)

| Scenarios / External costs | 1 | 2 | 3 | 4 | 5 | 6 |
|---|---|---|---|---|---|---|
| $ec_1$ | 0 | $1.4 \times 10^{-3}$ | $2.8 \times 10^{-3}$ | $4.3 \times 10^{-3}$ | $6.7 \times 10^{-3}$ | $7.1 \times 10^{-3}$ |
| $ec_2$ | 0 | $0.7 \times 10^{-3}$ | $1.4 \times 10^{-3}$ | $2.1 \times 10^{-3}$ | $2.8 \times 10^{-3}$ | $3.5 \times 10^{-3}$ |
| $ec_3$ | 0 | $0.005 \times 10^{-3}$ | $0.011 \times 10^{-3}$ | $0.017 \times 10^{-3}$ | $0.023 \times 10^{-3}$ | $0.028 \times 10^{-3}$ |

Table 6. Available fuel resources for each of PPs

|  | $\sigma_1$ (Million−lit) | $\sigma_2$ (Million−lit) | $\sigma_3$ (Million−m$^3$) |
|---|---|---|---|
| Available source of fuel for each PP per a year | 100 | 100 | 100 |

Table 7. Permissible amount for emission of pollutants - (m$^3$×10$^3$)

|  | $Z_1$ | $Z_2$ | $Z_3$ |
|---|---|---|---|
| Permitted amount of pollutant emission | 1,240 | 6,000 | 800,000 |

After introducing all parameters required creating a numerical example, then, we have studied this numerical example in the following using two meta-heuristic algorithms.

### 5.1. Computational results

To evaluate the numerical example, both algorithms have been coded in MATLAB 2016. Moreover, this numerical example has been solved in a laptop with Core i5 CPU and 8 GB of RAM. Before solving the problem, we have tuned the parameters of GA. To do this, we have applied the Response Surface Methodology (RSM) method using Minitab 17. For this purpose, four factors, i.e., parameters, such as the number of initial population, number of iteration, crossover rate, and mutation rate have been evaluated. Three levels have been assumed for each of these factors; these levels have been provided in Table (8). Note that we have originated PSO algorithm form Clerc and Kennedy (2002), and adapted it for this problem. For this reason, we use the parameter setting suggested by Clerc and Kennedy (2002) in the rest of this paper.

Table 8. Range of parameters settings for the proposed GA

| Levels<br>Parameters | 1 | 2 | 3 |
|---|---|---|---|
| Crossover rate | 0.7 | 0.8 | 0.9 |
| Mutation rate | 0.1 | 0.2 | 0.3 |
| Number of iteration | 800 | 1000 | 1200 |
| Number of initial population | 200 | 300 | 400 |

To tune the parameters of GA, we have studied one-way and two-way interactions. Fig. (3) shows the importance of these interactions with 90% of confidence level. Considering P-Values, $A$, $C$, $B \times B$, and $C \times D$ are identified as the important interactions. Thus, Eq. (21) may effectively tune and determine the best settings of factors.

>>Please insert Fig. (3) about here>>

$$y = 739533 - 6767A + 2094C + 4160B \times B + 2809C \times D \quad (21)$$

After optimizing Eq. (21), the best possible setting for the factors discussed earlier has been provided in Table (9). Based on Clerc and Kennedy (2002), we have also provided the optimal settings for PSO algorithm in this table.

Table 9. Parameter settings considered for the proposed algorithms

| The proposed GA | | The proposed PSO | |
|---|---|---|---|
| Initial population | 400 | Initial population | 400 |
| Max iteration | 1000 | Max iteration | 1000 |
| Crossover rate | 0.7 | $\varphi_1$ and $\varphi_2$ | 2.05 |
| Mutation rate | 0.2 | Personal learning coefficient | $C_1 = W \times \varphi_1$ |
| Parent selection | Tournament | Global learning coefficient | $C_2 = W \times \varphi_2$ |

After tuning the parameters of algorithms, each of 12 instances created in the numerical example, i.e., six instances for the scenarios of the external costs in the collusion market and six instances for the scenarios of the external costs in the competitive market, has been solved for 15 times using the proposed algorithms. Then, we have applied the $2 - Sample\ t$ test with 90% of confidence level in order to figure out which of these algorithms preforms better. To do so, we have determined the total production of PPs and created a sample consisting 15 data (due to solving each instance 15 times with each of algorithms) for each instance. Then, we have compared the data of a sample obtained by GA with its corresponding sample obtained by PSO algorithm. The results of $2 - Sample\ t$ test indicates that both algorithms perform identically for more than 90% of instances. In the rest of instances, PSO algorithm mostly outperforms GA. To sake for brevity, we have only illustrated the results obtained of the comparison between the proposed GA and PSO for the first scenario of external costs and in the collusion market. As illustrated in Fig. (4), there is no significant difference between the obtained results by the proposed GA and PSO. However, the standard deviation for the solutions provided by PSO algorithm is much smaller than the ones provided by GA. In other words, it can be argued that the solutions provided by PSO are more reliable compared to solutions provided by GA. Thus, we have used the solutions provided by the proposed PSO in the rest of paper.

As mentioned previously, this paper aims to compare the behavior of PPs in collusion and competitive markets and in the presence of external costs. Fig. (5) illustrates the generated electricity by PPs with respect to both markets and various scenarios of external costs. Due to being equipped with more efficient equipment, the second PP has generated more electricity than first and third ones in both markets and different scenarios (near to its maximum capacity). According to Fig. (5), amount of electricity generated by the second PP does not differ in both markets significantly. It indicates that efficiency of equipment would decrease the sensitivity of PPs to the type of markets and external costs. The first PP has generated more electricity in the collusion market. On the other hand, the third PP has generated more electricity in the competitive market. In addition, it shows that amount of electricity generated by the first and third PPs have slightly increased and decreased by the increment of external costs, respectively. This is because of the fact that the equipment of the first PP is more efficient than the third PP. It should be noted that PPs and markets have been indicated in Fig. (5) using some abbreviations. For instance, 'PCL1' and 'PCM1' refer to the first PP in collusion and competitive markets, respectively. The same kind of abbreviations have also been used in the rest of paper.

>>Please insert Fig. (4) about here>>

>>Please insert Fig. (5) about here>>

As illustrated in Fig. (6), profit acquired by all of PPs have decreased by the increment in the external costs. As the same as the production level, the second PP acquired more profit than the first and third ones. Considering that the first PP had generated more electricity in the collusion market, therefore, it has acquired more profit in this market compared to the competitive market. The third PP had generated more electricity in the competitive market, therefore, it has acquired more profit in this market compared to the collusion market. Although the first PP has generated more electricity in the sixth scenario of the collusion market than the competitive one, the increment in the external costs has decreased the margin profit of the first PP. In other words, generating more electricity in the presence of the external costs would not provide the same amount of profit as PPs could acquire in the absence of the external costs. The same trend can be argued for the third PP.

>>Please insert Fig. (6) about here>>

Fig. (7) shows the fuel resource consumption of PPs in both markets and in the presence of external costs. Based on Fig. (7), the consumption of gas has not differed by the increment in the external costs, and in both markets. This trend has occurred due to that fact that gas only emits a negligible amount of contaminants. On the other hand, the consumption of fuel-oil has decreased by the increment in the external costs, and in both markets. This trend has occurred due to that fact that fuel-oil emits a considerable amount of contaminants. Since production levels of all PPs are higher in the collusion market, and the available fuel resources of gas are limited, the second PP has used more gas-oil by the increment in the external costs. This trend has occurred because of the fact that equipment in this PP is more efficient compared to other PPs. Thus, this PP has consumed more gas-oil and less gas to let others use more gas. In this way, other PPs could use more gas, which is a greener fuel resource, in order to decrease the external costs. It should be noted that fuel resources, and markets have been indicated in Fig. (6) using some abbreviations. For instance, 'FOCL', 'GOCL', and 'GCM' refer to fuel-oil in the collusion market, gas-oil in the collusion market, and gas in the competitive market, respectively. Moreover, the consumption gas has been drawn based on the right axis, but the consumption of other fuel resources have been drawn based on the left axis.

>>Please insert Fig. (7) about here>>

Figs. (8) and (9) illustrate the emission of environmental pollutants, i.e., $CO_2$, $SO_2$, and $NO_x$ in both markets and in the presence of external costs. According to Fig. (8), emission of $SO_2$ in both markets has gradually decreased by the increment in external costs. Considering the fact that fuel-oil is the main source for the emission of $SO_2$, the decrement in the use of fuel-oil has decreased the emission of this contaminant. Likewise, the same trend has happened for $CO_2$. In the case of $NO_x$, both fuel-oil and gas-oil are recognized as the main sources of this contaminant. Due to the fact that there is not a specific trend for the consumption of these fuel resources when the external costs increase, therefore, the consumption of $NO_x$ does not also pursue a specific trend over various scenarios.

>>Please insert Fig. (8) about here>>

>>Please insert Fig. (9) about here>>

Eventually, Table (10) compares the total profit acquired by all PPs in both markets and in the various scenarios. According to Table (10), the total profit acquired by all PPs is slightly more in the collusion market compared to the competitive market. Moreover, it is clear that increment in the external costs has decreased the total profit acquired by all PPs in both markets.

Table 10. Profit share of PPs in the cases of collusion and competition (Billion USD)

| Mode of cooperation \ Scenarios | 1 | 2 | 3 | 4 | 5 | 6 |
|---|---|---|---|---|---|---|
| Total profit in the collusion market | 202.817 | 192.602 | 182.437 | 172.436 | 162.595 | 152.826 |
| Total profit in the competitive market | 202.613 | 192.404 | 182.246 | 172.250 | 162.434 | 152.645 |

Taking into account what has been discussed so far, it can be concluded that increment in the external costs would not significantly affect PPs that are equipped with efficient equipment in both markets. In other words, these PPs are not sensitive to the external costs. Therefore, the increment of the external costs could be a powerful incentive for PPs to equip themselves with more efficient equipment. In other words, governments can mitigate the environmental contaminants emitted by PPs in the presence of the external costs.

## 6. Conclusion

In this study, two mathematical formulations were proposed to study the behavior of PPs in the collusion and competitive markets and in the presence of external costs. To solve this model, two meta-heuristic algorithms, i.e., GA and PSO algorithm, were also proposed. The parameters of GA was tuned using the RSM method. Moreover, the parameters of PSO algorithm were also tuned based on Clerc and Kennedy (2002). Then, the performance of these two algorithms was compared using $2 - Sample\ t$ test. The results showed that the proposed PSO outperforms the proposed GA in terms of mean value and standard deviation of solutions. For this reason, this algorithm was used in the rest of the paper.

Afterward, a numerical example was created and introduced, which includes three different PPs regarding their equipment. These PPs have been evaluated in both collusion and competitive markets, and in the presence of the external costs. To determine how the presence of the external costs could affect the PPs, we defined six different scenarios for the external costs. Assessment of this numerical example showed that increment in the external costs affects PPs almost similarly in both markets. Whereas, PPs affect less in the collusion market compared to the competitive market. The computational results showed that increment in the external costs motivates PPs to use more gas-oil and less fuel-oil. Despite the lower cost of purchasing fuel-oil, PPs found it more efficient to use more gas-oil in the presence of the external costs. Increasing the external costs does not affect the consumption of gas in both markets due to its low contaminant emission rate.

The computational results also showed that not only profit of PPs would decrease in the presence of the external costs, but also the emission of environmental contaminants might decrease. According to the computational results, it can be concluded that increment in the external costs would not significantly affect PPs that are equipped with efficient equipment in both markets. In other words, these PPs are not sensitive to the external costs. Therefore, the increment of the external costs could be a powerful incentive for PPs to equip themselves with more efficient equipment. In other words, governments can mitigate the environmental contaminants emitted by PPs in the presence of the external costs.

## References

Bushnell, J. B., Mansur, E. T., & Saravia, C. (2007). Vertical Arrangements, Market Structure, and Competition an Analysis of Restructured US Electricity Markets: National Bureau of Economic Research.

Chen, H., Wong, K., Nguyen, D., & Chung, C. (2006). Analyzing oligopolistic electricity market using coevolutionary computation. Power Systems, IEEE Transactions on, 21(1), 143-152.

Clerc, M., & Kennedy, J. (2002). The particle swarm-explosion, stability, and convergence in a multidimensional complex space. IEEE transactions on Evolutionary Computation, 6(1), 58-73.

Cunningham, L. B., Baldick, R., & Baughman, M. L. (2002). An empirical study of applied game theory: Transmission constrained Cournot behavior. Power Systems, IEEE Transactions on, 17(1), 166-172.

El-Hawary, M. E. (1995). Electrical power systems: design and analysis (Vol. 2): John Wiley & Sons.

Ellersdorfer, I. (2005). A multi-regional two-stage Cournot model for analyzing competition in the German electricity market.

Deputy Director of Energy. (2013). Balance sheet of Iran's energy in 2011. Iran: Ministry of Power.

Evans, J., & Green, R. J. (2004). Why did British electricity prices fall after 1998?

Fabra, N., & Toro, J. (2005). Price wars and collusion in the Spanish electricity market. International Journal of Industrial Organization, 23(3), 155-181.

Fan, L., Norman, C. S., & Patt, A. G. (2012). Electricity capacity investment under risk aversion: A case study of coal, gas, and concentrated solar power. Energy Economics, 34(1), 54-61.

Georgakellos, D. A. (2012). Climate change external cost appraisal of electricity generation systems from a life cycle perspective: the case of Greece. Journal of Cleaner Production, 32, 124-140.

Grainger, J. J., & Stevenson, W. D. (1994). Power system analysis: McGraw-Hill.

Hobbs, B. F. (2001). Linear complementarity models of Nash-Cournot competition in bilateral and POOLCO power markets. Power Systems, IEEE Transactions on, 16(2), 194-202.

Deputy Director of Statistics and Information. (2012). Detailed stats of Iran's electricity industry, particularly for production in 2011. Iran: Mother company of Tavanir (in Persian).

Jaehnert, S., & Doorman, G. L. (2012). Assessing the benefits of regulating power market integration in Northern Europe. International Journal of Electrical Power & Energy Systems, 43(1), 70-79.

Jehle, G. A. (2001). Advanced microeconomic theory: Pearson Education India.

Kennedy, J., & Eberhart, R. (1995, Nov/Dec 1995). Particle swarm optimization. Paper presented at the Neural Networks, 1995. Proceedings., IEEE International Conference on.

Khoshakhlagh, R., Sharifi, A., & Parvan, H. (2013). Proposed method to calculate the social costs of generating electricity in power plants in short term (case study of Shahid Mohammad Montazeri and Isfahan Islamabad's power plants) Economic Modeling Research, 10, 77-97.

Ladjici, A., Tiguercha, A., & Boudour, M. (2014). Nash Equilibrium in a two-settlement electricity market using competitive coevolutionary algorithms. International Journal of Electrical Power & Energy Systems, 57, 148-155.

Lajevardi, H., & Mohades, N. (2010). Comparison of end costs of each Kwh electricity in different pricing ways. Seventh national conference on energy.

Li, G., Shi, J., & Qu, X. (2011). Modeling methods for GenCo bidding strategy optimization in the liberalized electricity spot market–A state-of-the-art review. Energy, 36(8), 4686-4700.

Li, T., & Shahidehpour, M. (2005). Strategic bidding of transmission-constrained GENCOs with incomplete information. Power Systems, IEEE Transactions on, 20(1), 437-447.

Lise, W., Linderhof, V., Kuik, O., Kemfert, C., Östling, R., & Heinzow, T. (2006). A game theoretic model of the Northwestern European electricity market—market power and the environment. Energy Policy, 34(15), 2123-2136.

Moitre, D. (2002). Nash equilibria in competitive electric energy markets. Electric power systems Research, 60(3), 153-160.

Roth, I. F., & Ambs, L. L. (2004). Incorporating externalities into a full cost approach to electric power generation life-cycle costing. Energy, 29(12), 2125-2144.

Shafie-Khah, M., Moghaddam, M. P., & Sheikh-El-Eslami, M. K. (2013). Development of a virtual power market model to investigate strategic and collusive behavior of market players. Energy Policy, 61, 717-728.

Stoft, S. (2002). Power system economics. Journal of Energy Literature, 8, 94-99.

Streimikiene, D., & Alisauskaite-Seskiene, I. (2014). External costs of electricity generation options in Lithuania. Renewable Energy, 64, 215-224.

Sweeting, A. (2007). Market power in the England and Wales wholesale electricity market 1995–2000*. The Economic Journal, 117(520), 654-685.

Vrhovcak, M. B., Tomsic, Z., & Debrecin, N. (2005). External costs of electricity production: case study Croatia. Energy Policy, 33(11), 1385-1395.

Weinzettel, J., Havránek, M., & Ščasný, M. (2012). A consumption-based indicator of the external costs of electricity. Ecological Indicators, 17, 68-76.

Wen, F., & David, A. (2001). Oligopoly Electricity Market Production under Incomplete Information. IEEE Power Engineering Review, 21(4), 58-61. doi: 10.1109/MPER.2001.4311314

Whitley, D. (1994). A genetic algorithm tutorial. Statistics and computing, 4(2), 65-85.

Willems, B. (2002). Modeling Cournot competition in an electricity market with transmission constraints. The Energy Journal, 95-125.

Willems, B., Rumiantseva, I., & Weigt, H. (2009). Cournot versus Supply Functions: What does the data tell us? Energy Economics, 31(1), 38-47.

Yucekaya, A. D., Valenzuela, J., & Dozier, G. (2009). Strategic bidding in electricity markets using particle swarm optimization. Electric power systems Research, 79(2), 335-345.

| 0.7490 | 0.6903 | 0.7005 | 0.6185 | 0.5418 | 0.6616 | 0.8170 | 0.6740 | 0.7490 |

Fig. 1. Solution representation of the proposed algorithms

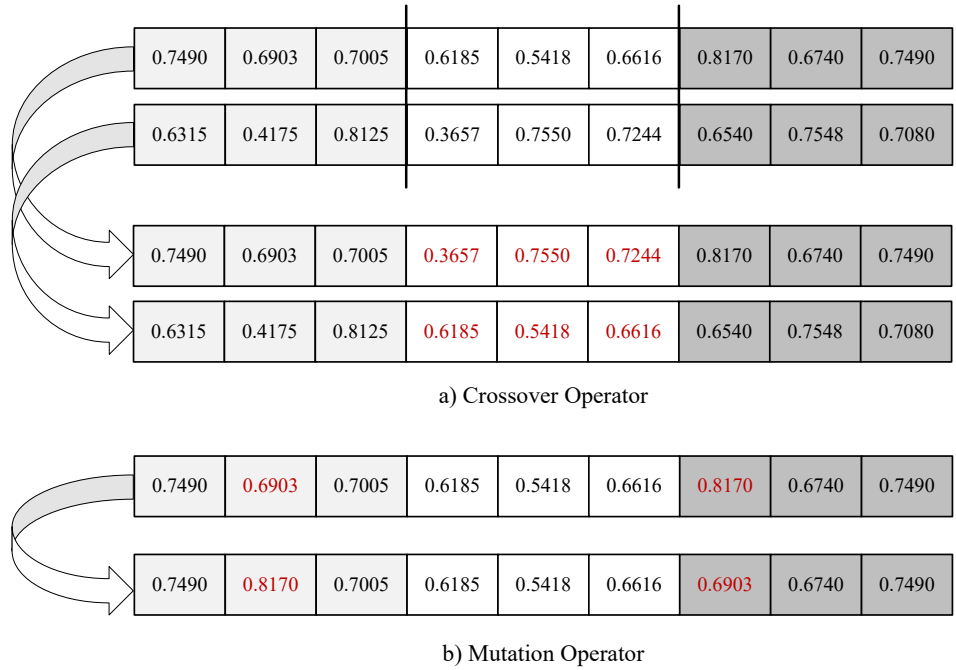

a) Crossover Operator

b) Mutation Operator

Fig. 2. The illustrative examples for the two-point crossover and swap mutation operators

```
Coded Coefficients

Term       Effect      Coef   SE Coef  T-Value  P-Value   VIF
Constant           7357479      2702  2723.35    0.000
A           -3869    -1934      1103    -1.75    0.084
B            3374     1687      1103     1.53    0.131   1.00
C            3758     1879      1103     1.70    0.093   1.00
D             814      407      1103     0.37    0.713   1.00
A*A          2361     1181      1910     0.62    0.539   1.00
B*B          8320     4160      1910     2.18    0.033   1.00
C*C         -4534    -2267      1910    -1.19    0.240   1.00
D*D          4555     2277      1910     1.19    0.237   1.00
A*B          1618      809      1351     0.60    0.551   1.00
A*C          2942     1471      1351     1.09    0.280   1.00
A*D         -4451    -2225      1351    -1.65    0.104   1.00
B*C           293      146      1351     0.11    0.914   1.00
B*D          3342     1671      1351     1.24    0.221   1.00
C*D          5619     2809      1351     2.08    0.041   1.00
```

Fig. 3. Effectiveness of factors and their interactions on objective function ($\alpha = 0.1$)

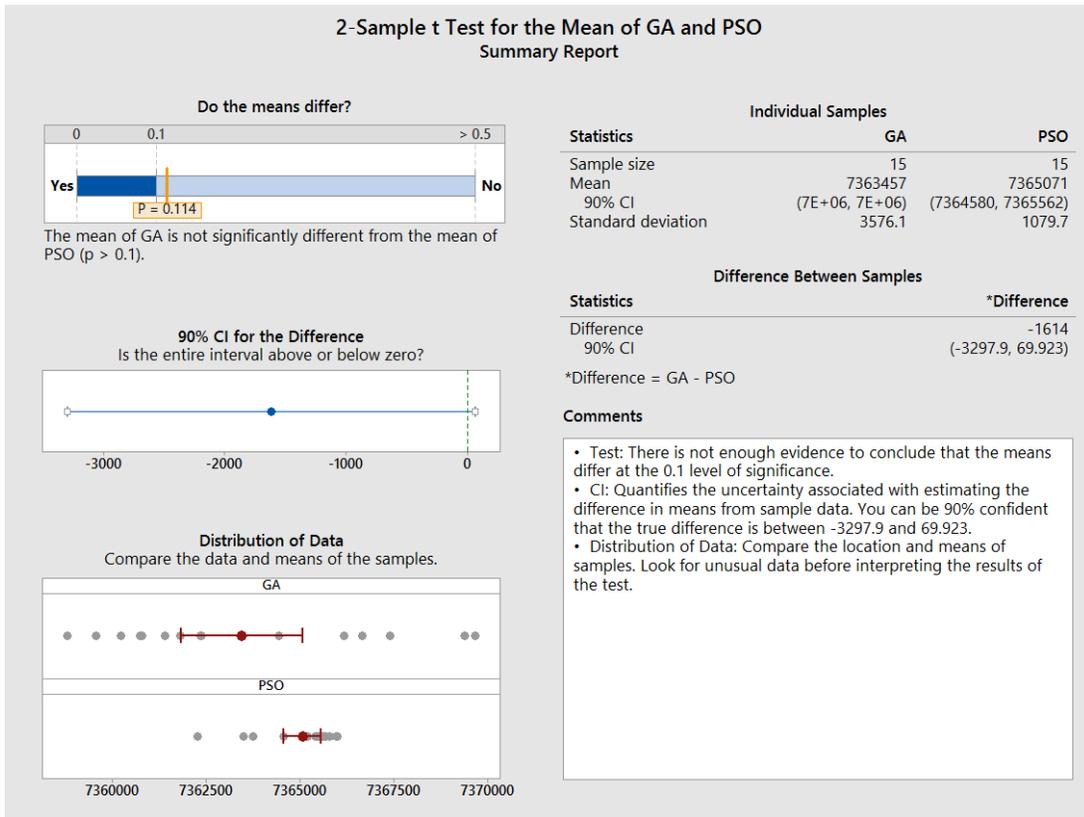

Fig. 4. Comparison of the results obtained by the proposed GA and PSO for the first scenario of external costs and in the collusion market

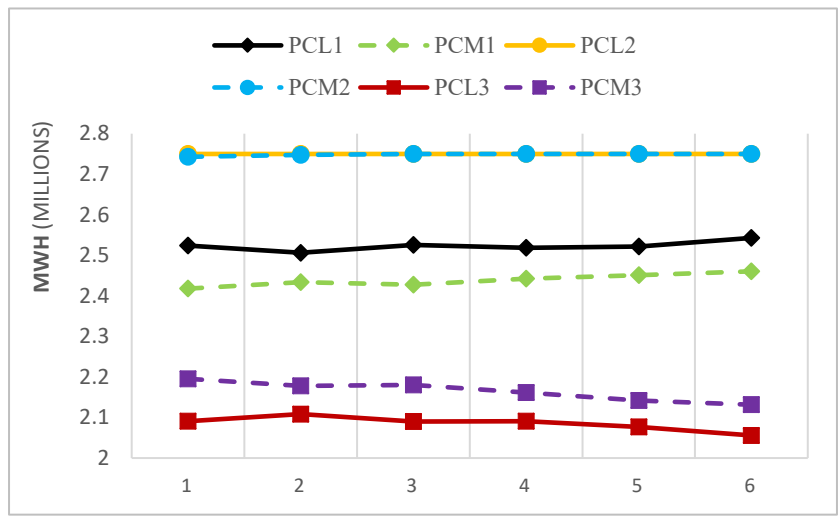

Fig. 5. Amount of electricity generated by PPs in the collusion and competitive markets and in the presence of external costs

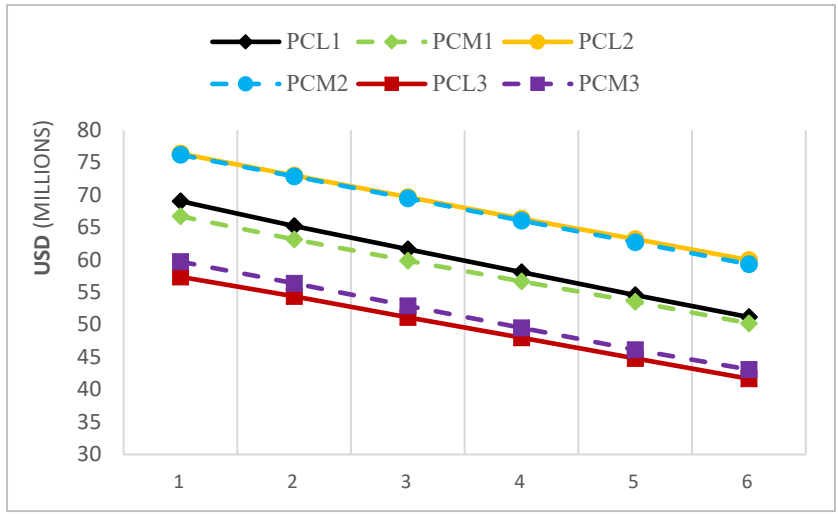

Fig. 6. Profit acquired by PPs in the collusion and competitive markets and in the presence of external costs

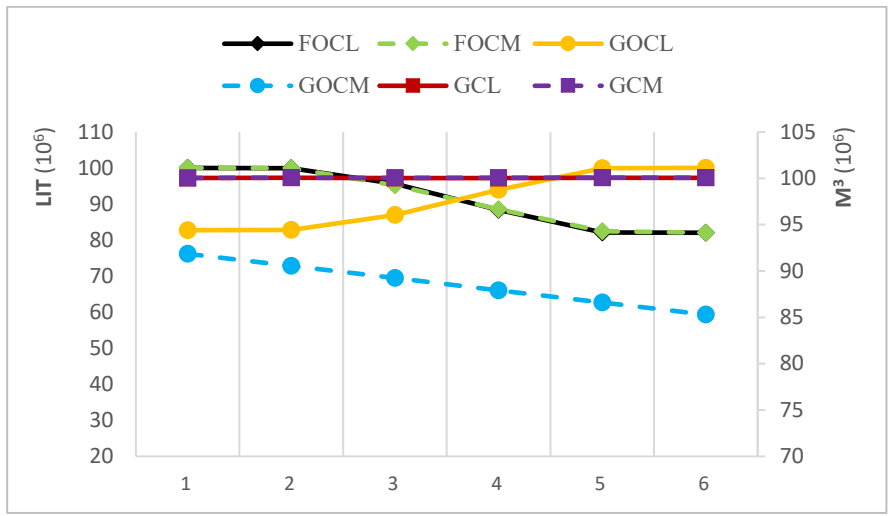

Fig. 7. The consumption of fuel resources by PPs in the collusion and competitive markets and in the presence of external costs

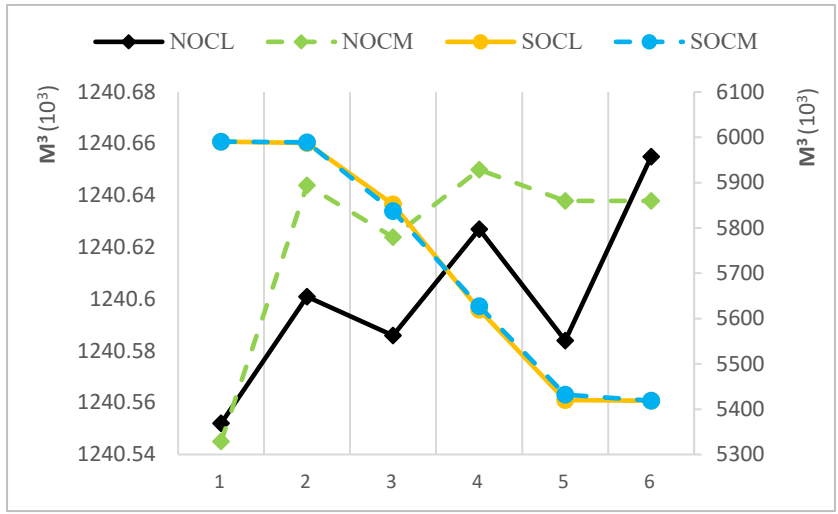

Fig. 8. Emission of $NO_x$ and $SO_2$ by PPs in the collusion and competitive markets and in the presence of external costs

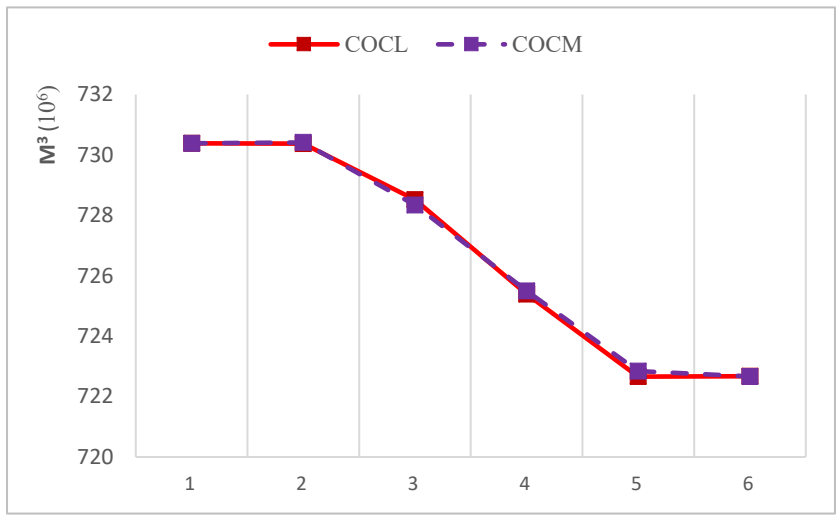

Fig. 9. Emission of $CO_2$ by PPs in the collusion and competitive markets and in the presence of external costs